\documentclass[12pt, reqno]{amsart}
\usepackage{amsmath, amsthm, amscd, amsfonts, amssymb, graphicx, color}
\usepackage[bookmarksnumbered, colorlinks, plainpages]{hyperref}
\usepackage{float}
\usepackage{graphicx,epstopdf}
\usepackage{bbm}
\usepackage{mathrsfs}
\usepackage{dsfont}
\usepackage{amssymb}
\usepackage{csquotes}
\usepackage[dvips]{epsfig}
\usepackage{latexsym}
\usepackage{exscale}
\usepackage[latin1]{inputenc}

\usepackage{tikz}
\usepackage{tikz-cd}
\usetikzlibrary{matrix,arrows,decorations.pathmorphing}

\hypersetup{colorlinks=true,linkcolor=red, anchorcolor=green, citecolor=cyan, urlcolor=red, filecolor=magenta, pdftoolbar=true}

\newtheorem{theorem}{Theorem}[section]
\newtheorem{lemma}[theorem]{Lemma}
\newtheorem{proposition}[theorem]{Proposition}
\newtheorem{corollary}[theorem]{Corollary}
\theoremstyle{definition}
\newtheorem{definition}[theorem]{Definition}
\newtheorem{example}[theorem]{Example}

\theoremstyle{remark}
\newtheorem{remark}[theorem]{Remark}
\numberwithin{equation}{section}

\newcommand{\norm}[1]{\lVert#1\rVert}

\newcommand{\R}{\mathbb{R}} %% Conjunto reales:        \R
\newcommand{\K}{\mathbb{K}} %% campo:     \K
\newcommand{\F}{\mathcal{F}}
\newcommand{\Lin}{\mathcal{L}}

\newcommand{\U}{\mathcal{U}}
\begin{document}

%%%%%%%%%%%%%%
\setcounter{page}{1}

\title{Schauder's Theorem and s-numbers}

\begin{center}
\author[A. G. AKSOY, D. A. THIONG ]{Asuman G\"{u}ven AKSOY, Daniel Akech Thiong}
\end{center}

\address{$^{*}$Department of Mathematics, Claremont McKenna College, 850 Columbia Avenue, Claremont, CA  91711, USA.}
\email{\textcolor[rgb]{0.00,0.00,0.84}{aaksoy@cmc.edu}}

\address{$^{1}$Department of Mathematics, Claremont Graduate University, 710 N. College Avenue, Claremont, CA  91711, USA..}
\email{\textcolor[rgb]{0.00,0.00,0.84}{daniel.akech@cgu.edu}}

\subjclass[2010]{Primary 47A16, 47B10; Secondary 47A68}

\keywords{s-numbers, approximation schemes, Schauder's theorem }

%\date{Received: xxxxxx; Revised: yyyyyy; Accepted: zzzzzz.
%\newline \indent $^{*}$Corresponding author}

\begin{abstract}
Motivated by the well known theorem of Schauder, we study the  relationship between various s-numbers of an operator $T$ and its adjoint $T^*$ between Banach spaces. \end{abstract} \maketitle

\section{Introduction}
%In the following we give a brief review of background, notation and terminology that will be relevant to this paper. 
Let $\mathcal{L}(X,Y)$ denote the normed vector space of all continuous operators from $X$ to $Y$, $X^*$ be the dual space of $X$, and $\mathcal{K}(X,Y)$ denote the collection of all compact operators from $X$ to $Y$. Denote by 
$T^{*} \in \mathcal{L}(Y^{*}, X^{*} )$ the adjoint operator of $T\in \Lin (X, Y)$. The well known theorem of Schauder states that $T \in \mathcal{K}(X,Y) \iff T^{*} \in \mathcal{K}(Y^{*},X^{*})$.  The proof of Schauder's theorem that uses Arzel$\grave{a}$-Ascoli Theorem is presented in most textbooks on functional analysis (see, e.g \cite {Rudin}).  A new and simple proof, which does not depend on Arzel$\grave{a}$-Ascoli can be found in \cite{Runde}.

Define  rank-$1$ operator $a\otimes y\in \mathcal{L}(X, Y)$ as $(a\otimes y) (x):=a(x)y \,\,\,\mbox{where}\,\,\, \\a\in X^*, \,\, y\in Y .$ An operator $T\in \mathcal{L}(X,Y)$ has finite rank if  $\ rank(T) := \dim \{Tx:\,\, x\in X\}$ is finite. Such an operator can be represented by $$ T=\sum_{k=1}^n a_k\otimes y_k\quad\mbox{with} \,\,a_1,\dots,a_n \in X^*\,\,\mbox{and}\,\, y_1,\dots, y_n\in Y.$$

 For two arbitrary normed spaces $X$ and $Y$, we define the collection of the finite-rank operators as follows: $$\F(X, Y) = \{A \in \mathcal{L}(X, Y): \text{rank} (A) \leq n -1 \},$$ which forms the smallest ideal of operators.

The concept of \textit{s-numbers}  $s_n(T)$  is introduced axiomatically in \cite{PieID}, and there are several possibilities  of  assigning  to every operator $T: X \to Y$ a certain sequence of numbers $\{s_n(T)\}$ such that 
$$ s_1(T) \geq s_2(T) \geq \dots \geq 0$$  which characterizes its degree of approximability or compactness of $T$. The main examples of s-numbers are approximation numbers, Gelfand numbers and Kolmogorov numbers. Among these  the largest  s-number is the  approximation numbers which we define below.  
%%%%%%%%%%%%%%%%%%%%
\begin{definition} 
The \textit{nth approximation number} $$\alpha_{n}(T) = \inf\{||T - A||: A \in \F(X, Y)\},\quad n=0,1,\dots$$
\end{definition}
$\alpha_{n}(T)$ provides a measure of how well T can be approximated by finite mappings whose range is at most n-dimensional.  Approximation numbers of an operator have the following properties \cite{PieID}:

\begin{enumerate}
\item $\alpha_0(T) = || T ||$
\item $\alpha_{n}(T) \geq \alpha_{n+1}(T)$ for all $n$
\item $\alpha_n( S+T) \leq \alpha_k(S) + \alpha_j(T)$ where $k+j=n$
\item $\alpha_n(\lambda T)=|\lambda| \alpha_n(T)$ for all $n$ and scalars $\lambda$

\item $|\alpha_n(S)-\alpha_n(T)| \leq ||S-T||$ for all $n$

\end{enumerate}

%%%%%%%%%%%%%%%%%%%%

%s-numbers are used to define quasi-Banach operator ideals $\mathcal{L}_{w}^{(\rho)}:=  \{ T\in \Lin(X,Y): \quad ( n^{\rho-1/w} s_n(T))\in \ell_w\}$ and their properties  and the relationship between s-numbers and the eigenvalue distribution are studied %by many. (See \cite{PieID} and \cite{Pieeig} and the references therein.)

\begin{definition} 
 We say that $T \in \mathcal{L}(X, Y)$ is of \textit{type $l^{p}$ }where $0<p< \infty$ if $(\alpha_n(T)) \in \ell^p$  or $(\alpha_n(T)) \in c_0  $ in case $p= \infty$.  For  $0<p< \infty$  in case  $\displaystyle\sum_{n}^{\infty}  (\alpha_{n}(T))^{p} < \infty$ and we denote such collection by $\ell^{p}(X,Y)$, which is again a linear subspace of $\mathcal{L}(X, Y)$ and it is the space of all linear operators of type $\ell^p$.

\end{definition}

s-numbers are used to define quasi-Banach operator ideals $$\mathcal{L}_{w}^{(\rho)}:=  \{ T\in \Lin(X,Y): \quad ( n^{\rho-1/w} s_n(T))\in \ell_w\}$$ and their properties  and the relationship between s-numbers and the eigenvalue distribution are studied by many. See for example  \cite{PieID} and \cite{Pieeig} and the references therein.

%%%%%%%%%%%%%%%%%
\section{Closure Results}
An \emph{operator ideal} $\U$ is a subclass of $\Lin$ such that the components $\U(E, F):=\U \cap \Lin(E, F)$ satisfy the following conditions:
%%%%%%%%%
\begin{enumerate}

\item [(i)] $I_{\K}\in \U$,

\item [(ii)] if $S_1, S_2\in \U(E, F)$ then $S_1+S_2\in \U(E, F)$,

\item [(iii)] if $T\in \Lin(E_0, E), S\in \U(E, F)$ and $R\in \Lin(F, F_0)$, then $RST\in \U(E_0, F_0)$.\\

\end{enumerate}
%%%%%%%%%%%%%

 Let $\U$ be an operator ideal, the map $\norm{\cdot}_\U: \U\to \R^+$ is called an \emph{$s$-norm on the operator ideal $\U$} with $0<s\le 1$ if the following are satisfied:
 
 \begin{enumerate}

\item [(i)] $\norm{I_{\K}}=1$,

\item [(ii)] $\norm{T_1+T_2}^s_\U\le \norm{T_1}^s_\U$ + $\norm{T_2}^s_{\U}$,

\item [(iii)]if $S\in \Lin (X_0, X)$, $T\in \U (X, Y)$ and $R\in \Lin (Y, Y_0)$ then $$\norm{RTS}_\U\le \norm{R} \norm{T}_\U \norm{S}$$
\end{enumerate}
 If each $\U(X, Y)$ is a Banach space with respect to the ideal $s$-norm $\norm{\cdot}_\U$, then $(\U, \norm{\cdot}_\U)$ is called an \emph{$s$-Banach operator ideal}.

\begin{definition} 
The \textit{dual ideal} $\U^{dual}$ consist of all operators $T$ such that $T^*$ belongs to the given ideal $\U$.   We say operator ideal is \textit{symmetric}  if $\U= \U^{dual}.$ Furthermore, we use the notation $\overline{{\U}}$, to denote the closed hull of the operator ideal  and the closure refers to the uniform topology on $\mathcal{L}(X,Y)$. 
\end{definition} 

More explicitly, a class of operators $\mathcal{A}(X,Y)\subset \mathcal{L}(X,Y)$ is called \textit{symmetric} if $T \in \mathcal{A}(X,Y) \implies T^{*} \in \mathcal{A}(Y^{*},X^{*})$. To be able to compare the degree of non-compactness of $T \in \mathcal{A}(X,Y)$ with that of $T^{*} \in \mathcal{A}(Y^{*},X^{*})$ requires $\mathcal{A}$ to be a \textit{symmetric ideal of operators}. 

The class $\mathcal{K}(X, Y)$ of compact operators between arbitrary Banach spaces $X$ and $Y$ is an example of a symmetric ideal of operators in $\mathcal{L}(X, Y)$.

Using the Principle of Local Reflexivity, Hutton (\cite{Hut} , Theorem 2.1) proved that  for $T \in \mathcal{K}(X, Y) $ implies  that $\alpha_{n}(T) = \alpha_{n}(T^{*})$ for all $n$.  However for non-compact operators $\alpha_{n}(T) \neq \alpha_{n}(T^{*})$ as shown in the following example:

\begin{example}  \cite{Hut}
Consider  $T = I: \ell_{1} \rightarrow c_{0}$  canonical injection  and $T^{*}: \ell_{1} \rightarrow \ell_{\infty}$  natural injection. Then, one has $\alpha_{n}(T) = 1$ for each $n$ and $\alpha_{n}(T^{*}) = \frac{1}{2}$.
\end{example}

For non-compact operator $T \in \mathcal{L}(X,Y)$, we do not have a lot of information about the relationship between $s_{n}(T)$ with $s_{n}(T^{*})$, however by imposing certain natural conditions on $X$ and $Y$ we were able to obtain relationship between  $s_{n}(T)$ with $s_{n}(T^{*})$ for certain s-numbers. 
 In this paper, we give a different proof of  Hutton's  result  using Kolmogorov diameters and for symmetrized  approximation numbers.  Additionally, we consider  operators which are not  compact  but compact with respect to certain approximation scheme, we call such operators $Q$-compact and prove a version of Schauder's theorem for $Q$-compact operators. In the case of non-compact operators, this answers the question of comparing the degree of compactness for $T$ and its adjoint $T^{*}$. 

\

The following theorem follows as a consequence of Hutton's  %(\cite{Hut} , Theorem 2.1)  
above mentioned theorem. %who proved that  $T \in \mathcal{K}(X, Y) \implies \alpha_{n}(T) = \alpha_{n}(T^{*})$ for all $n$.   

\begin{theorem}
The operator ideal $\overline{\F}$ is symmetric.
\end{theorem}

Pietsch has shown that the space of finite-rank linear operators is a dense subset of the space of all linear operators of type $\ell^{p}$ between Banach spaces (see, \cite{Piet72}, Prop. 8.2.5). This can be used to prove that every operator of type $\ell^{p} $ is relatively compact (see, \cite{Piet72}, Prop. 8.2.6) and hence compact since the notions coincide for Banach spaces. 

As a corollary, we have:

\begin{corollary}

 If $T \in  \overline{\mathcal{F}(X, Y)} = \ell^{p}(X,Y)$, then $\alpha_{n} (T) = \alpha_{n}(T^{*})$ for all $n$. 
 \end{corollary} 
 
 The corollary implies that for $0 < p \leq \infty$, $T \in \ell^{p}(X, Y) \iff T^{*} \in \ell^{p} (Y^{*}, X^{*})$, which shows that  $\ell^{p}(X,Y)$ is an example of a symmetric ideal of operators in $\mathcal{L}(X, Y)$. 

\section{ Hutton's Theorem Revisited}
In this section, we re-state a version of Hutton's theorem and give a  different proof which  uses the  basic theorems of functional analysis, together with  Principle of Local Reflexivity.  Lindenstrass and Rosenthal \cite{LR}  discovered a principle which shows that all Banach spaces are ``locally reflexive" or said in another way  every bidual $X^{**}$ is finitely representable  in the original space $X$. The following  is a stronger version of this property called \textit{Principle of Local Reflexivity} (PLR) due to Johson, Rosental and Zippin \cite{JRZ}:

\begin{definition}
Let $X$ be  a Banach space regarded  as a subspace of $X^{**}$, let $E$ and $F$  be finite dimensional subspaces of $X^{**}$ and $X^*$ respectively and  let $\epsilon >0$. Then there exist a one-to-one operator $T: E \to X$ such that
\begin{enumerate}
\item $T(x)= x$  for all $x\in X \cap E$
\item $f(Te)=e(f)$ for all $e\in E$ and $f\in F$
\item $||T|| ||T^{-1}|| < 1+\epsilon$.
\end{enumerate}
\end{definition}
 PLR is an effective tool in Banach space theory and  more recently   Oja and Silja  in \cite{Oja}  investigated versions of the  principle of local reflexivity for nets of subspaces of  a Banach space and  gave  some applications to some duality and lifting theorems.
 Next, we define  Kolmogoroff diameter of $T \in \mathcal{L}(X)$  and observe an alternate way of characterizing compact operators using Kolmogoroff diameter  of $T$.  
\begin{definition} [\cite{CS}, Prop. 2.2.2]
The \textit{nth -Kolmogoroff diameter} of $T \in \mathcal{L}(X)$ is defined by $$\delta_{n}(T) = \inf \{||Q_{G} T||: \dim G \leq n \}$$ where the infimum is over all subspaces $G \subset X$ and $Q_{G}$ denotes the canonical quotient map $Q_{G}: X \rightarrow X/G$. 

\end{definition} 

\begin{lemma} [Lemma 1 in \cite{Runde}]
Let $X$ be a Banach space and let $T \in \mathcal{L}(X)$. Then $T \in \mathcal{K}(X)$ if and only if, for each $\epsilon >0$, there is a finite-dimensional subspace $F_{\epsilon}$ of $X$ such that $||Q_{F_{\epsilon}}T|| < \epsilon$, where $Q_{F_{\epsilon}}: X \rightarrow X/F_{\epsilon}$. 

\end{lemma} 

In the following, we restate  Hutton's theorem and give a  different proof which  uses the  basic theorems of functional analysis, together with  PLR.

\begin{theorem} 
Let $T \in \mathcal{K}(X)$. Then $\alpha_{n} (T) = \alpha_{n} (T^{*})$. 

\end{theorem} 
\begin{proof}
Since one always has $\alpha_n(T^*)\leq \alpha_n(T)$, if we have $\alpha_n(T)\leq \alpha_n(T^{**})$, then $\alpha_n(T^{**})\leq \alpha_n(T^{*})$ would imply $\alpha_n(T)\leq \alpha_n(T^{*})$. 
Thus we must verify $\alpha_n(T)\leq \alpha_n(T^{**})$. 
To this end, suppose $T \in \mathcal{K}(X)$, by Schauder's theorem, $T^{*}$ and $T^{**}$ are compact. Let $\epsilon > 0$, then by definition, there exists $A \in \mathcal{F}(X^{**})$ such that $||T^{**} - A|| < \alpha_{n} (T^{**}) + \epsilon$. 

By Lemma 3.3,  there are finite-dimensional subspaces $E_{\epsilon}$ of $X^{**}$ and $F_{\epsilon}$ of $X^{*}$ such that $||Q_{E_{\epsilon}}T^{**}|| < \epsilon$, where $Q_{E_{\epsilon}}: X^{**} \rightarrow X^{**}/E_{\epsilon}$ and  $||Q_{F_{\epsilon}}T^{*}|| < \epsilon$, where $Q_{F_{\epsilon}}: X^{*} \rightarrow X^{*}/F_{\epsilon}$.

By the Principle of Local Reflexivity (PLR), there exists a one-to-one linear operator $S: E_{\epsilon} \rightarrow X$ such that $||S||||S^{-1}|| < 1 + \epsilon$, $y^{*}(Sx^{**}) = x^{**}(y^{*})$ for all $x^{**} \in E_{\epsilon}$ and all $y^{*} \in F_{\epsilon}$, and $S_{| E_{\epsilon} \cap X} = I$.

Let $J: X \to X^{**}$ be the canonical map. By the Hahn-Banach theorem, since $E_{\epsilon} $ is a subspace of $X^{**}$, $S: E_{\epsilon} \rightarrow X$ can be extended to a linear operator $\overline{S}: X^{**} \rightarrow X$. 

We now have $T \in \mathcal{L}(X)$ and $\overline{S}AJ \in \mathcal{L}(X)$ and $\text{rank }(\overline{S}AJ ) = \text{rank} (A) < n$, and therefore $$\alpha_{n} (T) \leq ||T - \overline{S}AJ ||.$$

To get an upper bound for $||T- \overline{S}AJ||$  we estimate $||Tx - \overline{S}AJ(x)||$ for $x \in B_{X}$ using an appropriate element $z_{j}$ of the covering of the set $T(B_{X})$.

Indeed, the compactness of $T$ implies that $T(B_{X})$ is relatively compact so that one can extract a finite-dimensional subset $Y_{\epsilon} \subset T(B_{X}) \subset X$ and let $z_{j} = Tx_{j}$ be the n elements forming a basis. 

Let $x \in B_{X}$. Then we have $$ \alpha_{n} (T) \leq |Tx - \overline{S}AJ(x)||  \leq ||Tx - z_{j} || + ||z_{j} - \overline{S}AJ(x)|| $$ $$\leq \epsilon + ||z_{j} - \overline{S}AJ(x)|| = \epsilon + ||\overline{ S} z_{j} -  \overline{S}AJ(x)||  \leq \epsilon + (1 + \epsilon) ||z_{j} - AJ(x) || <  \epsilon + (1 + \epsilon) (\alpha_{n} (T^{*}) + \epsilon),$$

since $$||z_{j} - AJ(x) || = ||Jz_{j} - AJ(x)|| \leq ||Jz_{j}  - JTx|| + ||JTx -  AJ(x)|| \leq \epsilon + ||JTx - AJx||$$

$$ = \epsilon + ||T^{**} Jx - AJx|| \leq ||T^{**} - A|| < \alpha_{n} (T^{*}) + \epsilon.$$

It follows that  $\alpha_{n} (T) \leq \alpha_{n} (T^{**})$, as  promised.

\end{proof} 
\begin{remark}
Since a  nuclear operator is compact for which a trace may be defined (nuclear operators on Hilbert spaces are called trace-class operators), it is natural to ask how nuclearity of $T$ and $T^*$ are related.
Recall that if  $T\in \mathcal{L}(X,Y)$ is a nuclear operator with the nuclear representation of $ T  = \displaystyle \sum_{n=1}^{\infty} \phi_n \otimes y_n$ then its adjoint  defined as  $T^*(\psi)= \displaystyle\sum_{n=1}^{\infty} \psi(y_n) \phi_n$ and  its nuclear norm defined as : $$||T||_{\mathcal{N}}= \inf\left \{ \displaystyle\sum_{n=1}^{\infty} ||\phi_n|| || y_n||: \,\,\,\, T(x) = \displaystyle\sum_{n=1}^{\infty}\phi_n(x) y_n \right\}$$   where the infimum is taken over all representations of $T$ of the form $T(x) = \displaystyle\sum_{n=1}^{\infty}\phi_n(x) y_n $  and $(\phi_n)$ and $(y_n)$ are bounded sequences in $X^*$ and $Y$ respectively satisfying  $\displaystyle\sum_{n=1}^{\infty} ||\phi_n|| || y_n|| < \infty$.  It is known that in case  $X^*$ has the approximation property and  if the operator $T\in \mathcal{L}(X,Y)$ has a nuclear adjoint, then $T$ is nuclear as well   and $||T||_{\mathcal{N}}= ||T^*||_{\mathcal{N}}$ (see Proposition $4.10$ in \cite{Ryan}). 
\end{remark}

%%%%%%%%%%%%%%%%%%%%%%%%%%
%%%%%%%%%%%%%%%%%%%%%%

\section{Compactness with Approximation schemes} 
Approximation schemes were introduced in Banach space theory by Butzer and Scherer in 1968 \cite{But} and independently  by Y. Brudnyi and N. Kruglyak under the name of ``approximation families'' in 1978 \cite{BK}. They were
popularized by Pietsch in his 1981 paper \cite{Pi}, for  later developments we refer the reader to  \cite{ Ak-Al, AA,  AL}. 

% The following characterization of compact operators is the motive for the introduction of the Kolmogrov numbers defined with respect to certain approximation scheme on $X$. 
 %%%%%%%%%%%%%%%%%%%%%%
 \iffalse
\begin{theorem} \cite{Piet65}

An operator $T \in \mathcal{L} (X, Y)$ between arbitrary Banach spaces $X$ and $Y$ is compact if and only if for every $\epsilon > 0$ there exists a finite-dimensional subspace $N_{\epsilon} \subset F$ with $T(B_{X}) \subset N_{\epsilon} + \epsilon B_{Y}$. 
\end{theorem} 
 
 Let $M$ be a bounded subset in normed space $X$. We define \textit{the nth Kolmogorov diameter} of $M$: \[\delta_{n}(M) = \inf \{ \epsilon > 0: \text{ there is a linear subspace } F \subset X \text{ with } \dim F < n \text { for which } M \subset \epsilon B_{X} + F \}\]
 
A bounded subset $M$ of normed space $X$ is relatively compact if and only if \[ \lim_{n \to \infty} \delta_{n}(M) = 0 \]

\fi 
\begin{definition}[Generalized Approximation Scheme]
Let $X$ be a Banach space.  For each $n\in \mathbb{N}$, let $Q_n=Q_n(X)$ be a family of subsets of $X$ satisfying the following conditions:
\begin{itemize}
\item[$(GA1)$] $\{0\}=Q_0\subset Q_1\subset \cdots \subset Q_n \subset \dots$.
\item[$(GA2)$]  $\lambda Q_n \subset Q_n$ for all $n \in N$ and all scalars $\lambda$.
\item[$(GA3)$]  $Q_n+Q_m \subseteq Q_{n+m}$ for every $n,m \in N$.
\end{itemize}
Then $Q(X)=(Q_n(X))_{n \in N}$ is called a \emph{generalized approximation scheme} on $X$.  We shall simply use $Q_n$ to denote $Q_n(X)$ if the context is clear.
\end{definition}
 We use here the term ``generalized'' because the elements of $Q_n$ may be subsets of $X$. Let us now give a few important examples of generalized approximation schemes.
\begin{example}

\
\begin{enumerate}
%\item The classical approximation schemes introduced in Pietsch in his seminal paper \cite{Pie}.
\item  $Q_n=$ the set of all at-most-$n$-dimensional subspaces of any given Banach space $X$.
\item Let $E$ be a Banach space and $X=L(E)$; let $Q_n=N_n(E)$, where $N_n(E)=$ the set of all $n$-nuclear maps on $E$ \cite {PieID}.
\item Let $a^k=(a_n)^{1+\frac{1}{k}},$ where $(a_n)$ is a nuclear exponent sequence. Then {$Q_n$}  on $X=L(E)$ can be defined as the set of all $\Lambda_\infty (a^k)$-nuclear maps on $E$ \cite{Dubinsky_Ram}.
\end{enumerate}
\end{example}
%We are now able to introduce $Q$-compact sets and operators:
\begin{definition}[Generalized Kolmogorov Number]
Let $B_X$ be the closed unit ball of $X$,  $Q= Q(X)=(Q_n(X))_{n \in N}$ be a \emph{generalized approximation scheme} on $X$,  and $D$ be a bounded subset of $X$.  Then the $n^{\text{th}}$ \emph{generalized Kolmogorov number} $\delta_n(D;Q)$ of $D$ with respect to $Q$ is defined by
\begin{equation}
\label{GenKolmogorovNumber}
\delta_n(D;Q)=\inf\{r>0:D \subset rB_X+A \text{ for some }A \in Q_n(X)\}.
\end{equation}
Assume that $Y$ is a Banach space and $T \in \mathcal{L}(Y,X)$. The $n^{\text{th}}$ Kolmogorov number $\delta_n(T;Q)$ of $T$ is defined as $\delta_n(T(B_Y);Q)$.
\end{definition}
It follows that $\delta_n(T;Q)$ forms a non-increasing sequence of non-negative numbers:
\begin{equation}
\|T\|=\delta_0(T;Q)\geq \delta_1(T;Q)\geq \cdots \geq \delta_n(T;Q)\geq 0.
\end{equation}
We are now able to introduce $Q$-compact sets and operators:
%\subsection{$Q$-compact sets and maps} 

\begin{definition}[$Q$-compact set]
Let $D$ be a bounded subset of $X$. We say that $D$ is $Q$-\emph{compact} if $\displaystyle\lim_n \delta_n(D;Q)=0$.
\end{definition}

\begin{definition}[$Q$-compact map] We say that $T\in L(Y,X)$ is a $Q$-\emph{compact map} if $\displaystyle\lim_n \delta_n(T;Q)=0$, i.e., $T(B_Y)$ is a $Q$-compact set.
\end{definition}
 There are examples of  $Q$-compact maps which are not compact, first such map involves projections $P: L_p[0,1] \to R_p $ where $R_p$ denotes the closure of  the span of the  space of Radamacher functions (see \cite{Asu} for details ), another example  is the weighted backward shift operator $B_w$ on $c_0( \mathbb{N})$ with $w=\{w_n\}$ not converging to $0$ is $Q$-compact but not compact.
 
 Our first objective here is whether or not Schauder's type of theorem is true for $Q$-compact maps. For this purpose we utilize other s-numbers such as Gelfand and  symmetrized approximation numbers of $T$, which we define below.
 %%%%%%%%%%%%%%%%%%%%%

 %%%%%%%%%%%%%%%%%%%%%
% First we give the definition of Gelfand numbers.
  \begin{definition} The nth Gelfand number $c_{n} (T)$ is defined as: \[c_{n}(T) = \inf \{ \epsilon > 0: ||Tx|| \leq \sup_{1 \leq i \leq k} | \langle x, a_{i} \rangle |+ \epsilon ||x||, \text {where }  a_{i} \in X^{*}, 1 \leq i \leq k \text { with }  k < n \} \]
\end{definition}
It follows that an operator $T$ is compact if and only if $c_{n}(T) \to 0$ as $n \to \infty$. 

  %%%%%%%%%%%%%%%%%%%%
  %The following characterization of compact operators due to Pietsch,  is closely connected with the definition of the so-called Gelfand numbers. 

%It follows that an operator $T$ is compact if and only if $c_{n}(T) \to 0$ as $n \to \infty$. 

%\begin{theorem} [\cite{CS}, Prop. 2.3.1]
%An operator $T \in \mathcal{L}(X, Y)$ between arbitrary Banach spaces $X$ and $Y$ is compact if and only if for $\epsilon > 0$ there are finitely many functionals $a_{i} \in X^{*}, 1 \leq i \leq n_{\epsilon}$, such that $$||Tx|| \leq \sup_{1 \leq i \leq n_{\epsilon} } |\langle x, a_{i} \rangle |+ \epsilon ||x|| \text{ for all } x \in X.$$
%\end{theorem} 

\

It is possible to compare various s-numbers such as $a_{n}(T), \delta_{n}(T), c_{n}(T)$ if one imposes some mild restrictions on $X$ and $Y$.

\begin{definition}

We say that a Banach space $X$ has \textit{the lifting property} if for every $T\in \mathcal{L}(X, Y/F)$ and every $\epsilon >0$ there exists an operator $S \in \mathcal{L}(X, Y)$ such that $||S|| \leq (1 + \epsilon) ||T||$ and $T = Q_{F}S$, where $F$ is a closed subspace of the Banach space $Y$ and $Q_{F}: Y \rightarrow Y/ F$ denotes the canonical projection. 

\end{definition} 

\begin{definition}

A Banach space $Y$ is said to have \textit{the extension property} if for each $T \in \mathcal{L}(M, Y)$ there exists an operator $S \in \mathcal{L}(X, Y)$ such that $T = SJ_{M}$ and $||T|| = ||S||$, where $M$ is a closed subspace of of an arbitrary Banach space $X$ and $J_{M}: M \rightarrow Y$ the canonical injection. 

\end{definition} 

\

Before defining symmetrized approximation numbers of $T$, we need to consider two universally important Banach spaces.

\

The Banach space $\ell_{1} (\Gamma)$ of \textit{summable number families} $\{ \lambda_{\gamma \in \Gamma} \}$ over an arbitrary index set $\Gamma$, whose elements $\{ \lambda_{\gamma \in \Gamma} \}$ are characterized by $\sum_{\gamma \in \Gamma} | \lambda_{\gamma }| < \infty$, has the metric lifting property. If $T$ is any map from a Banach space with metric lifting property to an arbitrary Banach space, then $a_{n}(T) = \delta_{n} (T)$ (cf. \cite{CS}, Prop. 2.2.3). It is known that every Banach space $X$ appears as a quotient space of an appropriate space $\ell_{1} (\Gamma)$ (for a proof of this, see \cite{CS}, p. 52). 

\

The Banach space $\ell_{\infty} (\Gamma)$ of \textit{bounded number families} $\{ \lambda_{\gamma \in \Gamma} \}$ over an arbitrary index set $\Gamma$ has the metric extension property. If $T$ is any map from an arbitrary Banach space into a Banach space with metric extension property, then $a_{n}(T) = c_{n} (T)$ (cf. \cite{CS}, Prop. 2.3.3). It is known that every Banach space $Y$ can be regarded as a subspace of an appropriate space $\ell_{\infty} (\Gamma)$ (for a proof of this, see \cite{CS}, p. 60).

\begin{remark}

If $T \in \mathcal{L}(X,Y)$, where X and Y are arbitrary Banach spaces with metric lifting and extension property, respectively, then $\delta_{n} (T) = a_{n}(T) = c_{n} (T)$. 
It is also known that if $T \in \mathcal{L}(X,Y)$, where X and Y are arbitrary Banach spaces, then $\delta_{n}(T^{*}) = c_{n}(T)$ (cf. \cite{CS}, Prop. 2.5.5).
Hence,  $\delta_{n}(T^{*}) = c_{n}(T) = a_{n}(T) = \delta_{n}(T)$, which we summarize as a theorem below.

\begin{theorem}
If $T \in \mathcal{L}(X,Y)$, where X and Y are arbitrary Banach spaces with metric lifting and extension property, respectively, then $\delta_{n}(T^{*}) = \delta_{n}(T)$ for all $n$. 
\end{theorem} 

\begin{remark}
Astala in \cite{As} proved that if $T \in \mathcal{L}(X,Y)$, where X and Y are arbitrary Banach spaces with metric lifting and extension property, respectively, then $\gamma (T) = \gamma (T^{*})$, where $\gamma (T)$ denotes the measure of non-compactness of $T$. In \cite{Ak-Al}, it is shown that $\lim_{n \to \infty} \delta_{n}(T) = \gamma (T)$. This relationship between Kolmogorov diameters and the measure of non-compactness together with theorem 4.10 provide an alternative proof for the result of Astala. 
\end{remark} 

If $T \in \mathcal{K}(X,Y)$, then it is known that $\delta_{n}(T) = c_{n}(T^{*})$ (cf. \cite{CS}, Prop. 2.5.6). If X and Y are Banach spaces with metric lifting and extension property, respectively, then we have $\delta_{n} (T) = a_{n}(T) = c_{n} (T)$. 
Thus, we have the following theorem. 
\begin{theorem}
If $T \in \mathcal{K}(X,Y)$, where X and Y are arbitrary Banach spaces with metric lifting and extension property, respectively, then $c_{n}(T^{*}) = c_{n}(T)$ for all $n$. 
\end{theorem} 

\end{remark} 

\begin{remark} [\cite{MS}]
It is known that if $X$ has the lifting property, then $X^{*}$ has the extension property. However, if $Y$ has the extension property, then $Y^{*}$ has the lifting property if and only if $Y$ is finite-dimensional. 
\end{remark} 

\begin{remark} If $X$ has the lifting property and $Y$ is finite-dimensional with the extension property, then by remark 4.10, $Y^{*}$ has the lifting property and $X^{*}$ has the extension property, so that by remark 4.9, we have $\delta_{n} (T^{*}) = a_{n}(T^{*}) = c_{n} (T^{*})$.
\end{remark}

For our needs, we choose the closed unit ball $B_{Z}$ of the Banach space $Z$ as an index set $\Gamma$.  Our proof of the Schauder's theorem for Q-compact operators will depend on the fact that $\ell_{1} (B_{Z})$ has the lifting property and $\ell_{\infty} (B_{Z})$ has the extension property.

 \begin{definition}
 The \textit{nth symmetrized approximation number} $\tau_{n}(T) $ for operator $T$ between arbitrary Banach spaces $X$ and $Y$ is defined as follows:
 
 $$\tau_{n}(T) = \delta_{n}(J_{Y} T), $$
 
where $J_{Y}: Y \to \ell_{\infty} (B_{Y^{*}})$ is an embedding map

 \end{definition} 
 
 \begin{remark}

Definition 4.15 is equivalent to $$\tau_{n}(T)  = a_{n}(J_{Y}TQ_{X})$$ as well as to  $$\tau_{n}(T)  =c_{n}(TQ_{X}), $$ where $Q_{X}: \ell_{1}(B_{X}) \rightarrow X$. 

\end{remark} 

 \begin{proposition} [Refined version of Schauder's theorem \cite{CS}, p. 84]
 An operator $T$ between arbitrary Banach spaces $X$ and $Y$ is compact if and only if $$\lim_{n\to \infty} \tau_{n}(T)  = 0 $$ and moreover, $$\tau_{n}(T) = \tau_{n}(T^{*}).$$

 \end{proposition}

 %%%%%%%%%%%%%%%%%%%%%%%%%%
 \iffalse
 
 \begin{Remark} [Exercise 6.1 in \cite{Bre}] 
 
 Let $$\lambda_{1} \geq \lambda_{2} \geq \cdots \lambda_{n} \geq \cdots \geq 0$$ be a non-increasing sequence of non-negative numbers and consider the diagonal operator $$D: \ell_{p} \rightarrow \ell_{p}: D (x_{1}, x_{2}, \cdots, x_{n}, \cdots) = (\lambda_{1} x_{1}, \lambda_{2} x_{2}, \cdots, \lambda_{n} x_{n}, \cdots) .$$ Then $D$ is compact if and only if $$\lim_{n \to \infty} \lambda_{n} = 0.$$
 \end{Remark} 
 
 \begin{proof}
 We give a proof here that relies on computing the Gelfand numbers , approximation numbers, Kolmogorov numbers, and symmetrized approximation numbers for identity operator and the relations among them. A computation done in [\cite{CS}, p.85] gives 
 $$\tau_{n}(D) = \alpha_{n}(D) = \lambda_{n}.$$ 
 
 Now by the refined version of Schauder's theorem, $$ D \text{ is compact } \iff \lim_{n \to \infty} \tau_{n} (D) = 0 \iff \lim_{n \to \infty} \lambda_{n} = 0.$$
 \end{proof} 
 \fi
 %%%%%%%%%%%%%%%%%%%%%%%%%%
 Motivated by this, we give the definition of Q-compact using the symmetrized approximation numbers. 
%It follows that an operator $T$ is compact if and only if $c_{n}(T) \to 0$ as $n \to \infty$. 
  
  %%%%%%%%%%%%%%%%%%%
  \begin{definition} 
  We say $T$ is Q- symmetric compact if and only if $$\lim_{n \to \infty} \tau_{n} (T, Q) = 0.$$
   \end{definition} 
   \begin{remark} [\cite{CS}, Prop. 2.5.4-6] 
   
   \
 
 \begin{enumerate}
 
 \item[a)]  From remark 4.16, we have $\tau_{n}(T, Q) = c_{n} (TQ_{X}, Q)$, where $Q_{X}: \ell_{1}(B_{X}) \rightarrow X$. 
 
 \item [b)] We will also abbreviate the canonical embedding $K_{\ell_{1}(B_{Y^{*}})}: \ell_{1} (B_{Y^{*}}) \rightarrow \ell_{\infty} (B_{Y^{*}})^{*} $ by $K$ so that $Q_{Y^{*}} = J^{*}_{Y}K$. 
 
 \item[ c)] Denote by $P_{0}: \ell_{\infty} (B_{X^{**}} ) \rightarrow \ell_{\infty} (B_{X})$ the operator which restricts any bounded function on $B_{X^{**}}$ to the subset $K_{X}(B_{X}) \subset B_{X^{**}} $ so that $Q^{*}_{X} = P_{0} J_{X^{*}}$. 
 
 \item [d)]  The relations (b) and (c) are crucial facts for the estimates of $\delta_{n}(T^{*}, Q^{*})$ and $c_{n}(T^{*}, Q^{*})$. In particular, we have $c_{n}(T^{*}, Q^{*}) \leq \delta_{n}(T, Q)$. 
\end{enumerate} 
 \end{remark}

  We now state and prove (adopting similar proof due to Pietsch and reproduced in \cite{CS}, Prop. 2.6) the following.  
  
  \begin{theorem} [ Schauder's theorem for Q-compact operators]
  An operator $T$ between arbitrary Banach spaces $X$ and $Y$ is Q- symmetric compact if and only if $$\lim_{n \to \infty} \tau_{n} (T, Q)  = 0$$ and moreover, $$\tau_{n} (T^{*}, Q^{*})  = \tau_{n} (T, Q), $$
  that is to say the degree of Q-compactness of $T$ and $T^*$ is the same in so far as it is measured by the symmetrized approximation numbers $\tau_{n}$. 
  \end{theorem} 
  
  \begin{proof}
  
 The first part is the definition. So it suffices to show $\tau_{n} (T^{*}, Q^{*})  = \tau_{n} (T, Q)$. By Remark 4. 19 (a) and (b) we have the following estimates: \[\tau_{n} (T^{*}, Q^{*})  = c_{n}(T^{*}Q_{Y^{*}}, Q^{*}) = c_{n}(T^{*}J_{Y}^{*} K, Q^{*}) \leq c_{n}((J_{Y}T)^{*}, Q^{*}) \leq \delta_{n} (J_{Y}T, Q) = t_{n}(T, Q)\] Conversely, we have by using Remark 4.19 (c) and (d): \[ t_{n}(T, Q) = c_{n}(TQ_{X}, Q) = \delta_{n}(TQ_{X})^{*}, Q^{*}) = \delta_{n}(Q_{X}^{*}T^{*}, Q^{*} ) \] \[= \delta_{n}(P_{0} J_{X^{*}} T^{*}, Q^{*}) \leq \delta_{n}(J_{X^{*}}T^{*}, Q^{*})  = t_{n} (T^{*}, Q^{*})\]

  \end{proof} 

Next we define approximation numbers with respect to a given scheme as follows:
 \begin{definition}
Given an approximation scheme $\{Q_n\}$  on $X$ and $ T \in \mathcal{L}(X)$, the n-th approximation number  $\alpha_n(T, Q)$ with respect to this approximation scheme is defined as:
$$ \alpha_n(T, Q)=\inf \{ ||T- B|| : \,\, B\in \mathcal{L}(X), \,
\,\, B(X) \in Q_n\}$$

\end{definition}

Let $X^*$ and $X^{**}$ be the dual and second dual of $X$.  Note that if we let $J: X \rightarrow X^{**}$ be the canonical injection and let $(X, Q_{n})$ be an approximation scheme, then $(X^{**}, J(Q_{n}))$ is an approximation scheme.

Let $\{Q_n\}$ and $\{Q^{**}_n\}:=  \{ J(Q_{n}) \}$  denote the subsets of $X$ and  $X^{**}$ respectively. 

\begin{definition}
We say $(X, Q_n)$ has the \textit{extended local reflexivity property}  (ELRP) if for each countable subset $C$ of   $X^{**}$ , for each $F\in Q_n^{**}$ for some $n$  and each $\epsilon>0$, there exists a continuous linear map 
$$P: \mbox{span} (F \cup C) \to X$$
such that 
\begin{enumerate}
\item $||P|| \leq 1+\epsilon$
\item $P \restriction_{C\cap X}= I (Identity)$
\end{enumerate}
\end{definition}
Note that ELRP is an analogue of local reflexivity principle which is possessed by all Banach spaces. %  given in \cite{Hut}

\begin{theorem} 
Suppose $(X, Q_n)$  has ELRP and $T\in \mathcal{L}(X)$ has  separable range. Then for each $n$ we have  $\alpha_n(T, Q)= \alpha_n(T^*, Q^*)$.
\end{theorem}

\begin{proof}
Since one always have $\alpha_n(T^*, Q^*)\leq \alpha_n(T, Q)$ we only need to verify $\alpha_n(T, Q)\leq \alpha_n(T^{**}, Q^{**})$. Let $J: X \to X^{**}$ be the canonical map and $U_X$ be the unit ball of $X.$ Given $\epsilon> 0$, choose $B\in \mathcal{L}(X^{**})$ such that $B(X^{**}) \in Q^{**}_n$ and
$$ ||B-T^{**}|| < \epsilon + \alpha_n(T^{**}, Q^{**}_n).$$
Let $\{z_j\}$ be a countable dense set in $T(X)$, thus $Tx_j=z_j$ where $x_j\in X$. Consider  the set $$ K= \mbox{span} \{ (JTx_j)_1^{\infty} \cup B(X^{**})\}$$  applying  ELRP  of $X$ we obtain a map
$$ P : K \to X\,\,\mbox{such that} \,\,\,||P|| \leq 1+\epsilon \,\, \,\, \mbox{and}\,\,\, P \restriction_{(JTx_j)_1^{\infty} \cap X} =I$$

%%%%%%%%%%%%%%%
For $x\in U_X$, consider 
\begin{flalign*}
||Tx-PBJx|| &\leq ||Tx-z_j|| +||z_j- PBJx|| &\\
                   &\leq  \epsilon + || PJTx_j-PBJx||&\\
                   &\leq   \epsilon +(1+\epsilon) ||JTx_j- BJx|| &\\
                   &\leq  \epsilon +(1+\epsilon)[  ||JTx_j-JTx|| + || JTx-BJx||] &\\
                   &\leq  \epsilon +(1+\epsilon) [\alpha_n(T^{**}, Q^{**}_n) +2 \epsilon ]
  \end{flalign*} 
  
  and thus $$ \alpha_n(T, Q) \leq \alpha_n(T^{**}, Q^{**}_n).$$

%%%%%%%%%%%%%%%

\end{proof}

%Hutton proves the above theorem using approximation numbers $\alpha_{n}(T)$, thus we can refer to the operators on $\overline{{\F}}$ as \textit{approximable}. In \cite{Hut} it is shown that, for each compact operator $T \in \mathcal{L}(X,Y) $ where $X$ and $Y$ are Banach spaces, the equality of the approximation numbers  $\alpha_{n }(T)$ and  $\alpha_{n} (T^{*}) $ for each $n$; considering the identity map of $\ell_1$ into $c_0$ she shows the falsity of the above result for non-compact operators. More precisely,  if $T=  I: \ell_1 \to c_0$ natural injection, then $\alpha_{n}(T)=1$ for each $n$ but  for $T^{*}=I  :\ell_1 \to \ell_{\infty}$ , $\alpha_{n}(T^{*}) = \frac{1}{2}$  for each $n$.

%%%%%%%%%%%%%%%%%%%%%%

%%%%%%%%%%%%%%%%%%%%%%%%%%%%%%
 \bibliographystyle{amsplain}

\end{document}